# A NEW GENERALIZATION OF GENERALIZED HYPERGEOMETRIC FUNCTIONS

ARJUN K. RATHIE

In this paper a natural generalization of the familiar $H$-function of Fox namely the $I$-function is proposed. Convergence conditions, various series representations, elementary properties and special cases for the $I$-function have also been given.

## 1. Introduction.

The well known $H$-function of one variable, defined by Fox (1961) and studied by Braaksma (1964), contains as particular cases most of the special functions of applied Mathematics, but it does not contain some of the important functions such as the Riemann Zeta functions, polylogarithms etc. In Section 2, we mention a few examples of functions, which are not included in the $H$-function and hence suggesting the form of a generalization of the $H$ function. In Section 3, we shall define the generalization of the $H$-function, namely the "I-function" which contains the polylogarithms, the exact partition of the Gaussian model from statistical mechanics and functions useful in testing hypothesis from statistics as special cases. The existence of the contour integrals in these special cases were not studied by earlier workers. In Section 4, convergence conditions for this function have been derived. In Section 5, we list special cases of our function giving relations with the functions available in the literature including $H$ and $G$-functions. In Section 6, we give series representations in several cases

---

and in particular behaviour of the function for small values. In Section 7, we mention a few important properties.

## 2. Functions which are not particular cases of the $H$-function.

In statistical mechanics, one encounters Mellin-Barnes integrals of the type

$$(2.1) \qquad (2\pi i)^{-1} \int_L \phi(s) z^s \, ds,$$

where

$$(2.2) \qquad \phi(s) = \frac{\prod_{j=1}^{m} \Gamma^{B_j}(b_j - \beta_j s) \prod_{j=1}^{n} \Gamma^{A_j}(1 - a_j + \alpha_j s)}{\prod_{j=m+1}^{q} \Gamma^{B_j}(1 - b_j + \beta_j s) \prod_{j=n+1}^{p} \Gamma^{A_j}(a_j - \alpha_j s)}$$

where $A_j$, $j = 1, \ldots, p$, and $B_j$, $j = 1, \ldots, q$ are not, in general positive integers. Clearly for non-integral values of $A_j$ and /or $B_j$, (2.1) is not expressible as a $H$-function.

For examples:

(a) The Mellin-Barnes integral representation in the case of the free energy of a Gaussian model of phase transition, see Joycee (1972), in statistical mechanics as given by Inayat-Hussain (1987)

$$(2.3) \qquad \beta F(d; \varepsilon) = -\frac{1}{4\pi^{d/2}(1+\varepsilon)^2} (2\pi i)^{-1} \cdot$$

$$\cdot \int_{-i\infty}^{i\infty} \frac{\Gamma(-s)\Gamma^2(1+s)\Gamma^d(3/2+s)}{\Gamma^{1+d}(2s)} (-(1+\varepsilon)^{-2})^s \, ds,$$

where $i = \sqrt{-1}$ and $d \ (> 0)$ can take non-integer values.

(b) The Mellin-Barnes integral representation in the case of the Feynman integral $g$, for non-integer $m$, as expressed in the following form, see Inayat-Hussain (1987):

$$(2.4) \quad g(\tau, n, \mu, m, z) =$$

$$= \frac{K_{a-1} 2^{-m-2} \Gamma(m+1) \Gamma(1 + \tfrac{1}{2}\mu) B(\tfrac{1}{2}, \tfrac{1}{2} + \tfrac{1}{2}\mu)}{\pi \Gamma(\tau) \Gamma(\tau - \tfrac{1}{2}\mu)} \cdot$$

$$\cdot (2\pi i)^{-1} \int_{-i\infty}^{i\infty} \frac{\Gamma(-s)\Gamma(\tau + s)\Gamma(\tau - \tfrac{1}{2}\mu + s)}{(n+s)^{1+m} \Gamma(1 + \tfrac{1}{2}\mu + s)} (-z)^s \, ds.$$

A further example of a function which is not a special case of the $H$-function is the polylogarithm of complex order. For details, see Marichev (1983).

(c) Consider a $p$-variate random sample of size $N$ from the normal distribution with mean $\mu$ and covariance matrix $\Sigma$. The density function of $\lambda = L^{2/N}$, where $L$ is the likelihood ratio criterion for testing the hypothesis, see Anderson (1984)

$$(2.5) \qquad H_0 : \Sigma = \Sigma_0 \quad \text{and} \quad \mu = \mu_0.$$

($\Sigma_0$ is a given positive definite matrix and $\mu_0$ a given vector) is given by Nagarsenker and Pillai (1973, 74).

$$(2.6) \qquad f(\lambda) = K(p,n)\lambda^{N/2-1}(2\pi)^{p/2} \sum_{r=0}^{\infty} B_r \cdot$$

$$\cdot (2\pi i)^{-1} \int_{c-i\infty}^{c+i\infty} \lambda^{-t} t^{-\nu-r} \, dt$$

$$(2.7) \qquad = k(p,n)\lambda^{N/2-1}(2\pi)^{p/2} \sum_{r=0}^{\infty} B_r \cdot$$

$$\cdot (2\pi i)^{-1} \int_{c-i\infty}^{c+i\infty} \lambda^{-t} \frac{\Gamma^{\nu+r}(t)}{\Gamma^{\nu+r}(t+1)} \, dt,$$

where $\nu = p(p+3)/4$. Note that $\nu$ is a positive integer only for certain values of $p$.

A few other statistical problems, which give the density function in a similar form as mentioned above, are discussed by Korin (1968) and Nagarsenkar & Pillai (1973, 1974).

Thus the above examples in physics and statistics suggest a new generalization of the $H$-function which is given in the next section.

## 3. The $I$-function.

The $I$-function will be defined and represented by the following Mellin-Barnes type contour integral

$$(3.1) \qquad I_{p\,q}^{m\,n}\left[z \left| \begin{array}{c} (a_1, \alpha_1, A_1), \ldots, (a_p, \alpha_p, A_p) \\ (b_1, \beta_1, B_1), \ldots, (b_q, \beta_q, B_q) \end{array} \right. \right] =$$

$$= (2\pi i)^{-1} \int_L \phi(s) z^s \, ds,$$

where $\phi(s)$ is given by (2.2).
Also:

(i) $z \neq 0$;
(ii) $i = \sqrt{-1}$;
(iii) $m, n, p, q$ are integers satisfying $0 \leq m \leq q, 0 \leq n \leq p$;
(iv) $L$ is a suitable contour in the complex plane;
(v) an empty product is to be interpreted as unity;
(vi) $\alpha_j, j = 1, \ldots, p; \beta_j, j = 1, \ldots q; A_j, j = 1, \ldots, p;$ and $B_j$, $j = 1, \ldots, q$ are positive numbers;
(vii) $a_j, j = 1, \ldots, p$ and $b_j, j = 1, \ldots, q$ are complex numbers such that no singularity of $\Gamma^{B_j}(b_j - \beta_j s), j = 1, \ldots, m$ coincides with any singularity of $\Gamma^{A_j}(1 - a_j + \alpha_j s), j = 1, \ldots, n$. In general these singularities are not poles.

There are three different contours $L$ of integration.

(a) $L$ goes from $\sigma - i\infty$ to $\sigma + i\infty$, ($\sigma$ real) so that all the singularities of $\Gamma^{B_j}(b_j - \beta_j s), j = 1, \ldots, m$, lie to the right of $L$, and all the singularities of $\Gamma^{A_j}(1 - a_j + \alpha_j s), j = 1, \ldots, n$, lie to the left of $L$.
(b) $L$ is a loop beginning and ending at $+\infty$ and encircling all the singularities of $\Gamma^{B_j}(b_j - \beta_j s), j = 1, \ldots, m$, once in the clock-wise direction, but none of the singularities of $\Gamma^{A_j}(1 - a_j + \alpha_j s), j = 1, \ldots, n$.
(c) $L$ is a loop beginning and ending at $-\infty$ and encircling all the singularities of $\Gamma^{A_j}(1 - a_j + \alpha_j s), j = 1, \ldots, n$, once in the anti-clockwise direction, but none of $\Gamma^{B_j}(b_j - \beta_j s), j = 1, \ldots, m$.

We notice that the condition (a) can be interpreted as a particular case of the condition (c). Moreover in case (b), the sign of the function is changed with respect to the case (c). Hence when more than one of the definitions make sense, they lead to the same result.

In short, (3.1) will be denoted by

$$I_{pq}^{mn}\left[ z \left| \begin{array}{l} {}_1(a_j, \alpha_j, A_j)_p \\ {}_1(b_j, \beta_j, B_j)_q \end{array} \right. \right].$$

For $A_j$ (and / or $B_j$) not an integer, the poles of the gamma functions of the numerator in (2.2) are converted to branch points. The branch cuts can be chosen so that the path of integration can be distorted for each of the three contours $L$ mentioned above as long as there is no coincidence of poles from any $\Gamma(b_j - \beta_j s)$ and $\Gamma(1 - a_j + \alpha_j s)$ pair.

The convergence conditions of (3.1) are derived in the next section.

## 4. The convergence conditions.

The sufficient conditions for convergence of (3.1) are given in the following theorems:

**Theorem 1.** *The integral* (3.1) *for L defined by* (a), *converges when* $|\arg z| < \Delta \pi/2$, *if* $\Delta > 0$, *where*

$$(4.1) \quad \Delta = \sum_{j=1}^{m} B_j \beta_j - \sum_{j=m+1}^{q} B_j \beta_j + \sum_{j=1}^{n} A_j \alpha_j - \sum_{j=n+1}^{p} A_j \alpha_j.$$

*If* $|\arg z| = \Delta \pi/2$, $\Delta \geq 0$, *the integral* (3.1) *converges absolutely when*

(i) $\mu = 0$ *if* $\nabla > 1$, *where*

$$(4.2) \quad \mu = \sum_{j=1}^{q} B_j \beta_j - \sum_{j=1}^{p} A_j \alpha_j,$$

*and*

$$(4.3) \quad \nabla = \sum_{j=1}^{p} A_j [\operatorname{Re}(a_j) - 1/2] - \sum_{j=1}^{q} B_j [\operatorname{Re}(b_j) - 1/2];$$

(ii) $\mu \neq 0$, *if with* $s = \sigma + it$, $\sigma$ *and t real,* $\sigma$ *is chosen so that for* $|t| \to \infty$, *we have* $\nabla + \sigma \mu > 1$.

*Proof.* We apply Luke (1969)

$$(4.4) \quad |\Gamma(x + iy)| \sim (2\pi)^{1/2} e^{-\pi |y|/2} |y|^{x-1/2}, \quad (|y| \to \infty)$$

for $x, y$ real, to various gamma functions of the integrand of (3.1) with $s = \sigma + it$. Thus we obtain, for $|t| \to \infty$,

$$(4.5) \quad |\phi(s) z^s| \sim C |t|^{-\nabla - \sigma \mu} \exp[-t \arg z - \pi |t| \Delta/2],$$

where $C$ is independent of $t$. Hence Theorem 1 follows.

**Theorem 2.** *The integral* (3.1) *for L defined by* (b) *converges for* $q \geq 1$ *and either* $\mu > 0$ *or* $\mu = 0$ *if* $|z| < \nu$, *where*

$$\nu = \prod_{j=1}^{q} \beta_j^{B_j \beta_j} \bigg/ \prod_{j=1}^{p} \alpha_j^{A_j \alpha_j}. \tag{4.6}$$

*Proof.* We apply Dixon and Ferrar (1936)

$$|\Gamma(s+\alpha)| \sim (2\pi)^{1/2} e^{-\sigma - t\theta} r^{\sigma + \alpha - 1/2} \tag{4.7}$$

with $s = \sigma + it = re^{i\theta}$, $|\theta| < \pi - \varepsilon$, for any fixed small $\varepsilon$ and for any fixed real constant $\alpha$; to various gamma functions of the integrand of (3.1) with $s = \sigma + it$, $\sigma > 0$. Simplifying, we arrive at the following expression:

$$|\phi(s) z^s| \sim C |z e^\mu / \nu|^\sigma \sigma^{-\sigma \mu}, \tag{4.8}$$

where $C$ is independent of $\sigma$. The Theorem 2 immediately follows if we take $\sigma \to \infty$.

In a similar fashion, the following theorem can be established.

**Theorem 3.** *The integral* (3.1) *for the contour L defined by* (c) *converges when* $p \geq 1$ *and either* $\mu < 0$ *or* $\mu = 0$ *if* $|z| > \nu$.

## 5. Special cases.

It is easy to verify that

$$g(\tau, n, \mu, m; z) = \tag{5.1}$$

$$= \frac{K_{a-1} 2^{-m-2} \Gamma(m+1) \Gamma(1+\tfrac{1}{2}\mu) B(\tfrac{1}{2}, \tfrac{1}{2}+\tfrac{1}{2}\mu)}{\pi \Gamma(\tau) \Gamma(\tau - \tfrac{1}{2}\mu)} \cdot$$

$$\cdot I_{33}^{13} \left[ -z \,\bigg|\, \begin{array}{l} (1-\tau, 1, 1), (1-\tau+\tfrac{1}{2}\mu, 1, 1), (1-n, 1, 1+m) \\ (0, 1, 1), (-\tfrac{1}{2}\mu, 1, 1), (-n, 1, 1+m) \end{array} \right];$$

$$\beta F(d; \varepsilon) = -\frac{1}{4\pi^{d/2}(1+\varepsilon)^2} \cdot \tag{5.2}$$

$$\cdot I_{22}^{12} \left[ -(1+\varepsilon)^{-2} \,\bigg|\, \begin{array}{l} (0, 1, 2), (-1/2, 1, d) \\ (0, 1, 1), (-1, 1, 1+d) \end{array} \right];$$

(5.3) $$f(\lambda) = K(p,n)\lambda^{(N/2)-1}(2\pi)^{1/2p} \sum_{r=0}^{\infty} B_r \cdot$$
$$\cdot I_{1\,1}^{1\,0}\left[\lambda \left| \begin{array}{c} (1, 1, \nu+r) \\ (0, 1, \nu+r) \end{array}\right.\right];$$

(5.4) $\overline{H}_{pq}^{mn}[z] =$
$$= \overline{I}_{pq}^{mn}\left[z \left| \begin{array}{c} (a_1, \alpha_1, A_1), \ldots, (a_n, \alpha_n, A_n), \\ (b_1, \beta_1, 1), \ldots, (b_m, \beta_m, 1), \\ , (a_{n+1}, \alpha_{n+1}, 1), \ldots, (a_p, \alpha_p, 1) \\ , (b_{m+1}, \beta_{m+1}, B_{m+1}), \ldots, (b_q, \beta_q, B_q) \end{array}\right.\right].$$

Also $\bar{H}$ has been very recently introduced in the literature, see Inayat-Hussain (1987).

(5.5) $$H_{pq}^{mn}\left[z \left| \begin{array}{c} {}_1(a_j, \alpha_j)_p \\ {}_1(b_j, \beta_j)_q \end{array}\right.\right] = I_{pq}^{mn}\left[z \left| \begin{array}{c} {}_1(a_j, \alpha_j, 1)_p \\ {}_1(b_j, \beta_j, 1)_q \end{array}\right.\right]$$

Also $H_{pq}^{mn}$ is the $H$-function, see Fox (1961), Braaksma (1964), Mathai and Saxena (1978).

(5.6) $$G_{pq}^{mn}\left[z \left| \begin{array}{c} {}_1(a_j)_p \\ {}_1(b_j)_q \end{array}\right.\right] = I_{pq}^{mn}\left[z \left| \begin{array}{c} {}_1(a_j, 1, 1)_p \\ {}_1(b_j, 1, 1)_q \end{array}\right.\right]$$

Also $G_{pq}^{mn}$ is the $G$-function, see Luke (1969).

## 6. Series expansions.

This section deals with series expansions for the $I$-function in several cases. The following procedures, if applicable, may be employed depending on the numerical values of the parameters of the integrand in (3.1).

**Procedure 1.** The $I$-function (3.1) for $L$ defined by (a), can be expressed in a series form suitable for numerical computation under the following set of conditions:

(i) $a_j$, $j = 1, \ldots, p$, and $b_j$, $j = 1, \ldots, q$, are real;
(ii) $\mu = 0$.

From (3.1), we have

$$\ln \phi(s) = \tag{6.1}$$

$$= \sum_{j=1}^{m} B_j \ln \Gamma(b_j - \beta_j s) - \sum_{j=m+1}^{q} B_j \ln \Gamma(1 - b_j + \beta_j s) +$$

$$+ \sum_{j=1}^{n} A_j \ln \Gamma(1 - a_j + \alpha_j s) - \sum_{j=n+1}^{p} A_j \ln \Gamma(a_j - \alpha_j s).$$

Employing Luke (1969)

$$\ln \Gamma(z + a) = (z + a - \frac{1}{2})\ln z - z + \frac{1}{2}ln(2\pi) + \tag{6.2}$$

$$+ \sum_{k=1}^{r} \frac{(-1)^{k+1} B_{k+1}(a)}{k(k+1)2^k} + 0(z^{-r-1})$$

for $|\arg z| \leq \pi - \varepsilon$, $\varepsilon > 0$ and simplifying by using condition (ii), we find that $\phi(s)$ is of the form

$$\phi(s) = AB^s s^{-\nabla}\left[1 + \sum_{j=1}^{\infty} c_j s^{-j}\right] \tag{6.3}$$

where $c_j$ involve Bernoulli polynomials, $B = (-1)^{-\frac{\Delta}{2}} \nu^{-1}$ and

$$A = (-1)^{\sum_{j=1}^{m} B_j(b_j - \frac{1}{2}) - \sum_{j=n+1}^{p} A_j(a_j - \frac{1}{2})} (2\pi)^{\frac{1}{2}(\sum_{j=1}^{m} B_j - \sum_{j+m+1}^{q} B_j + \sum_{j=1}^{n} A_j - \sum_{j=n+1}^{p} A_j)} \tag{6.4}$$

$$\cdot \prod_{j=1}^{q} \beta_j^{B_j(b_j - \frac{1}{2})} \prod_{j=1}^{p} \alpha_j^{A_j(\frac{1}{2} - a_j)}$$

Explicit expression for $c_j$, in general, is too complicated and is therefore not given here. Thus (3.1) and (6.3) yield

$$I(z) = A(2\pi i)^{-1} \int_L (Bz)^s s^{-\nabla}(1 + \sum_{j=1}^{\infty} c_j s^{-j}) ds. \tag{6.5}$$

The Residue Theorem is applicable to the right hand side of (6.5) if $\nabla$ is a positive integer. In this case, we have

$$(6.6) \qquad I(z) = A \sum_{j=0}^{\infty} c_j [ln(Bz)]^{j+\nabla-1}/(j+\nabla-1)!$$

where $c_0 = 1$ and $Bz > 0$. If $\nabla$ is not a positive integer, then we can express $I(z)$ in a series of Beta functions for $0 < z < 1$.

**Procedure 2.** The $I$-function when $L$ is defined as in (a), $A'_j s$, $j = 1, \ldots, n$, and $B'_j s$, $j = 1, \ldots, m$, are positive integers, can be expressed in a series form suitable for numerical computation for the case when no singularity of $\prod_{j=m+1}^{q} \Gamma^{B_j}(1-b_j+\beta_j s)$ coincides with any poles of $\prod_{j=1}^{n} \Gamma^{A_j}(1-a_j+\alpha_j s)$ and no singularity of $\prod_{j=n+1}^{p} \Gamma^{A_j}(a_j-\alpha_j s)$ coincides with any poles of $\prod_{j=1}^{m} \Gamma^{B_j}(b_j-\beta_j s)$.
In this case

$$I(z) = \begin{cases} \text{sum of the residues of } \phi(s)z^s \text{ at the poles of} \\ \quad \prod_{j=1}^{m} \Gamma^{B_j}(b_j - \beta_j s), \ |z| < 1 \\ \text{sum of the residues of } \phi(s)z^s \text{ at the poles of} \\ \quad \prod_{j=1}^{n} \Gamma^{B_j}(1 - a_j + \alpha_j s), \ |z| > 1 \end{cases}$$

In particular, we shall find the series representation and behaviour for small values for the function

$$(6.7) \qquad \overline{I}(z) = \overline{I}_{pq}^{mn}\left[ z \left| \begin{array}{c} {}_1(a_j, \alpha_j, A_j)_{n}, {}_{n+1}(a_j, \alpha_j, A_j)_p \\ {}_1(b_j, \beta_j, 1)_{m}, {}_{m+1}(b_j, \beta_j, B_j)_q \end{array} \right. \right] =$$

$$= (2\pi i)^{-1} \int_L \frac{\prod_{j=1}^{m} \Gamma(b_j - \beta_j s) \prod_{j=1}^{n} \Gamma^{A_j}(1 - a_j + \alpha_j s)}{\prod_{j=m+1}^{q} \Gamma^{B_j}(1 - b_j + \beta_j s) \prod_{j=n+1}^{p} \Gamma^{A_j}(a_j - \alpha_j s)} z^s \, ds.$$

When the poles of $\Gamma(b_j - \beta_j s)$, $j = 1, \ldots, m$, are simple, the integral (6.7) can be evaluated with the help of the Residue Theorem to give

$$(6.8) \qquad \overline{I}(z) = \sum_{r=0}^{\infty} \sum_{h=1}^{m} \frac{\prod_{j=1}^{n} \Gamma^{A_j}(1 - a_j + \alpha_j \frac{b_h+r}{\beta_h})}{\prod_{j=n+1}^{p} \Gamma^{A_j}(a_j - \alpha_j \frac{b_h+r}{\beta_h})}.$$

$$\cdot \frac{\prod_{\substack{j=1 \\ \neq h}}^{m} \Gamma(b_j - \beta_j \frac{b_h+r}{\beta_h})(-1)^r z^{(b_h+r)/\beta_h}}{\prod_{j=m+1}^{q} \Gamma^{B_j}(1 - b_j + \beta_j \frac{b_h+r}{\beta_h}) r! \beta_h}$$

for $|z| < 1$.

From (6.8), it is clear that

(6.9) $\quad \overline{I}(z) \sim z^c, \quad \text{where} \quad c = \min_{i \leq j \leq m} [Re(b_j/\beta_j)]$

for small values of $z$.

On suitable specialising the parameters in (6.7), it is now easy to see that

(6.10) $\quad g(\tau, n, \mu, m; z) =$

$$= \frac{K_{a-1} 2^{-m-2} \Gamma(m+1) B(\frac{1}{2}, \frac{1}{2} + \frac{1}{2}\mu)}{\pi} \cdot$$

$$\cdot \sum_{r=0}^{\infty} \frac{(\tau - \frac{1}{2}\mu)_r (\tau)_r (n+r)^{-(1+m)}}{(1 + \frac{1}{2}\mu)_r} \frac{z^r}{r!}$$

for $|z| < 1$, and

(6.11) $\quad \beta F(d; \varepsilon) = -(1+\varepsilon)^{-2} 2^{-d-2} \sum_{r=0}^{\infty} \frac{(1)_r [(3/2)_r]^d}{[(2)_r]^{1+d} (1+\varepsilon)^{2r}}.$

Equations (6.10) and (6.11) are alternate forms which may be suitable for numerical computations of the functions in (5.1) and (5.2) respectively.

**Remark.** It may be remarked that the assumptions in this section exclude the polygarithms as well as many cases of (2.1).

**Procedure 3.** The $I$-function when $L$ is as defined in (b) and $B_j$, $j = 1, \ldots, m$, are positive integers, can be expressed in a series form suitable for numerical computation for the following two cases:

*Case 1.* No singularity of $\prod_{j=n+1}^{p} \Gamma^{A_j}(a_j - \alpha_j s)$ coincides with any poles of $\prod_{j=1}^{m} \Gamma^{B_j}(b_j - \beta_j s)$. In this case, $I(z) = $ sum of the residues of $\phi(s) z^s$ at the poles of $\prod_{j=1}^{m} \Gamma^{B_j}(b_j - \beta_j s)$.

*Case 2.* $A_j$, $j = n+1, \ldots, p$, are positive integers such that some of the poles of $\prod_{j=n+1}^{p} \Gamma^{A_j}(a_j - \alpha_j s)$ coincides with the poles of $\prod_{j=1}^{m} \Gamma^{B_j}(b_j - \beta_j s)$. In this case $I(z) = $ sum of the residues of $\phi(s)z^s$ at the poles of $\prod_{j=1}^{m} \Gamma^{B_j}(b_j - \beta_j s) \Big/ \prod_{j=n+1}^{p} \Gamma^{A_j}(a_j - \alpha_j s)$.

**Procedure 4.** The $I$-function when $L$ is a defined in (c) and $A_j$, $j = 1, \ldots, n$, are positive integers, can be expressed in a series form suitable for numerical computation for the following two cases:

*Case 1.* No singularity of $\prod_{j=m+1}^{q} \Gamma^{B_j}(1 - b_j + \beta_j s)$ coincides with the poles of $\prod_{j=1}^{n} \Gamma^{A_j}(1 - a_j + \alpha_j s)$. In this case, $I(z) = $ sum of the residues of $\phi(s)z^s$ at the poles of $\prod_{j=1}^{n} \Gamma^{A_j}(1 - a_j + \alpha_j s)$.

*Case 2.* $B_j$, $j = m+1, \ldots, p$ are positive integers such that some of the poles of $\prod_{j=m+1}^{q} \Gamma^{B_j}(1 - b_j + \beta_j s)$ coincide with the poles of $\prod_{j=1}^{n} \Gamma^{A_j}(1 - a_j + \alpha_j s)$. In this case $I(z) = $ sum of the residues of $\phi(s)z^s$ at the poles of $\prod_{j=1}^{n} \Gamma^{A_j}(1 - a_j + \alpha_j s) \Big/ \prod_{j=m+1}^{q} \Gamma^{B_j}(1 - b_j + \beta_j s)$.

The results in general, adopting the Procedure 2,3 or 4, will be obtained in terms of Gamma, Psi and generalized Zeta functions. These methods have been applied to various problems involving the derivation of the exact distribution of the likelihood ratio criterion, see Anderson (1984). Among the first papers in this direction, see Mathai and Rathie (1970, 1971). It may not be possible to apply Cauchy's residue theorem in other possible cases which have not been discussed here. In such cases one may have to employ other methods to evaluate the integral in (3.1).

**Remarks.**
1) The sufficient convergence condition for the function defined by (6.7) have recently been obtained by Buschman & Srivastava (1990).
2) For another generalization of the $H$-function, see Rathie (1989, 1994).

## 7. Simple properties.

The properties given below are immediate consequence of the definition (3.1) and hence they are given here without proof:

(i) The $I$-function is symmetric in the set of triplets

$$\{(a_1, \alpha_1, A_1), \ldots, (a_n, \alpha_n, A_n)\},$$
$$\{(a_{n+1}, \alpha_{n+1}, A_{n+1}), \ldots, (a_p, \alpha_p, A_p)\},$$
$$\{(b_1, \beta_1, B_1), \ldots, (b_m, \beta_m, B_m)\} \quad \text{and}$$
$$\{(b_{m+1}, \beta_{m+1}, B_{m+1}), \ldots, (b_q, \beta_q, B_q)\} \quad \text{respectively.}$$

(ii) If one of the triplets $(a_j, \alpha_j, A_j)$, $j = 1, \ldots, n$, is equal to one of the triplets $(b_j, \beta_j, B_j)$, $j = m+1, \ldots, q$. (Or one of the triplets $(a_j, \alpha_j, A_j)$, $j = n+1, \ldots, p$ is equal to one of the triplets $(b_j, \beta_j, B_j)$, $j = 1, \ldots, m$, then the $I$-function reduces to another $I$-function of the lower order. For example:

$$(7.1) \quad I_{p\,q}^{m\,n}\left[z \,\middle|\, \begin{array}{l} _1(a_j, \alpha_j, A_j)_p \\ _1(b_j, \beta_j, B_j)_{q-1}, (a_1, \alpha_1, A_1) \end{array}\right] =$$

$$= I_{p-1\,q-1}^{m\,n-1}\left[z \,\middle|\, \begin{array}{l} _2(a_j, \alpha_j, A_j)_p \\ _1(b_j, \beta_j, B_j)_{q-1} \end{array}\right]$$

provided that $p \geq n \geq 1$ and $q \geq m+1$.

$$(7.2) \quad I_{p\,q}^{m\,n}\left[z \,\middle|\, \begin{array}{l} _1(a_j, \alpha_j, A_j)_{p-1}, (b_1, \beta_1, B_1) \\ _1(b_j, \beta_j, B_j)_q, \end{array}\right] =$$

$$= I_{p-1\,q-1}^{m-1\,n}\left[z \,\middle|\, \begin{array}{l} _1(a_j, \alpha_j, A_j)_{p-1} \\ _2(b_j, \beta_j, B_j)_q \end{array}\right]$$

provided that $q \geq m \geq 1$ and $p \geq n+1$.

(iii)

$$(7.3) \quad z^\sigma I_{p\,q}^{m\,n}\left[z \,\middle|\, \begin{array}{l} _1(a_j, \alpha_j, A_j)_p \\ _1(b_j, \beta_j, B_j)_q \end{array}\right] = I_{p\,q}^{m\,n}\left[z \,\middle|\, \begin{array}{l} _1(a_j + \sigma\alpha_j, \alpha_j, A_j)_p \\ _1(b_j + \sigma\beta_j, \beta_j, B_j)_q \end{array}\right]$$

(iv)

$$(7.4) \quad I_{p\,q}^{m\,n}\left[z \,\Big|\, \begin{matrix}{}_1(a_j, \alpha_j, A_j)_p \\ {}_1(b_j, \beta_j, B_j)_q\end{matrix}\right] = c I_{p\,q}^{m\,n}\left[z^c \,\Big|\, \begin{matrix}{}_1(a_j, c\alpha_j, A_j)_p \\ {}_1(b_j, c\beta_j, B_j)_q\end{matrix}\right]$$

for $c > 0$.

(v)

$$(7.5) \quad I_{p\,q}^{m\,n}\left[z \,\Big|\, \begin{matrix}{}_1(a_j, \alpha_j, A_j)_p \\ {}_1(b_j, \beta_j, B_j)_q\end{matrix}\right] = I_{q\,p}^{n\,m}\left[z^{-1} \,\Big|\, \begin{matrix}{}_1(1-b_j, \beta_j, B_j)_q \\ {}_1(1-a_j, \alpha_j, A_j)_p\end{matrix}\right]$$

(vi) The $I$-function is, in general, a many valued function of $z$ with a branch point at $z = 0$.

A further study of the $I$-function will form the subject matter of a subsequent paper.

**Acknowledgement.** The author is grateful to the referee for his useful suggestions.

*Department of Mathematics,*
*Dungar College (M.D.S. University),*
*Bikaner - 334001 Rajasthan (INDIA)*